\documentclass[12pt,leqno,twoside]{amsart}
\usepackage{amssymb}
\usepackage{latexsym}
\usepackage{verbatim}
\usepackage{epsfig}

\newtheorem{theorem}{Theorem}[section]
\newtheorem{corollary}[theorem]{Corollary}

\newtheorem{proposition}[theorem]{Proposition}
\theoremstyle{definition}
\newtheorem{definition}[theorem]{Definition}
\theoremstyle{remark}
\newtheorem{remark}[theorem]{\sc Remark}
\theoremstyle{remark}

\theoremstyle{remark}
\newtheorem{example}[theorem]{\sc Example}
\theoremstyle{remark}
\newtheorem{note}[theorem]{\sc Note}
\theoremstyle{remark}

\theoremstyle{remark}

\renewcommand{\Box}{\square}    %\diamond
\renewcommand{\Bbb}{\mathbb}
\newcommand{\cal}{\mathcal}

\newcommand{\rk}{\mathop{{\rm{rk}}}\nolimits}

\newcommand{\h}{{\rm{ht}}}

\newcommand{\Sing}{\mathop{{\rm{Sing}}}\nolimits}
\renewcommand{\th}{{\rm{th}}}
\newcommand{\id}{{\rm{id}}}

\newcommand{\grad}{\mathop{\rm{grad}}\nolimits}
\newcommand{\ity}{{\infty}}

\newcommand{\fin}{\hspace*{\fill}$\Box$}

\newcommand{\bC}{{\Bbb C}}

\newcommand{\bP}{{\Bbb P}}

\newcommand{\bZ}{{\Bbb Z}}
\newcommand{\bX}{{\Bbb X}}
\newcommand{\bQ}{{\Bbb Q}}

\newcommand{\bH}{{\Bbb H}}

\newcommand{\bL}{{\Bbb L}}

\newcommand{\cS}{{\cal S}}
\newcommand{\cW}{{\cal W}}

\newcommand{\cH}{{\cal H}}

\def\growarrow#1{
  \setbox1=\hbox{ $\scriptstyle #1$\ }

\mathop{\smash{\hbox to \wd1{\rightarrowfill}}
          \vphantom\rightarrow}\limits^{#1}}

\begin{document}
%\date{}

\title[Homotopy groups of complements and non-isolated singularities]{Homotopy groups of complements and non-isolated singularities} 

\author{Anatoly Libgober}
\address{A.L.: Dept. of Mathematics, University of Illinois at Chicago, 851 S. Morgan Str. Chicago, Ill, 60607,
 USA.}   
\email{libgober@math.uic.edu}

\author{Mihai Tib\u ar}  
\address{M.T.:  U.F.R. Math\' ematiques, Unit\' e Mixte de Recherche 8524 CNRS,
Universit\'e de Lille 1, \  59655 Villeneuve d'Ascq, France.}
\email{tibar@agat.univ-lille1.fr}

\thanks{This work was started during the visit by the first author at the University of Lille in spring 2000. He wants to express his appreciation to the University of Lille for support and hospitality and also to thank NSF for support 
during work on this project.}



%\date{\today}

%\begin{abstract}
%\end{abstract}

\maketitle

\section{Introduction}

It has been known for some time that the topology of non-isolated singularities, at least in some cases, has something to do with the ``position'' of the singularities (cf. \cite[p.164]{St}, \cite{Dimca1}, \cite{Dimca2}). The starting 
point of this work was an attempt to clarify this relationship.
We shall consider two situations in which non-isolated singularities arise: 
\par  (1) {\it polynomials in $\bC^{n+1}$ with a single atypical value},
\par (2) {\it germs of analytic functions with non-isolated 
singularities.} 
\smallskip
\par It turns out that in the first case, say for a polynomial
$f : \bC^{n+1}\to \bC$ having $0$ as the only atypical value, the first 
non-trivial (in the appropriate sense, cf. section \ref{nonisolated})
homology group of a fiber $F_t = \{ f=t\}$, for generic $t$, is related to the first non-trivial 
homotopy group of the complement of $F_0 \cap H$ in $H$, for a generic
linear section $H$ of appropriate dimension.
The study of the first non-trivial homotopy groups of complements 
$\bC^{n+1}\setminus V$, 
was started in \cite{Li}. If 
$V$ does not have singularities at infinity and  
the dimension of the singular locus of $V$ is equal to $k$, 
the first possibly non-trivial higher homotopy group 
is the group $\pi_{n-k}({\bC}^{n+1}\setminus V)$ (since $\pi_1 ({\bC}^{n+1}\setminus V)=\bZ$
and  $\pi_i ({\bC}^{n+1}\setminus V)=0$ for $i <n-k$; cf. \cite{Li}). 
Although the information on higher homotopy groups of ${\bC}^{n+1}\setminus V$
depends heavily on the information about homotopy groups of spheres, 
the first non-trivial homotopy group has algebro-geometric meaning
and depends on the local type and the position of singularities 
of $V$ (cf. \cite{Li}, \cite{Li-2}). This identification 
allows one to relate in a direct way the homology of a smoothing 
to the position of singularities of a generic hyperplane section of the atypical fiber.
\par In the case (2) we obtain a 
similar relation.
More precisely, we relate the first non-trivial homology group of the
Milnor fiber of an analytic function $f$ (with non-isolated  singularity
at the origin) to the corresponding homotopy group of the complement to the zero set $f^{-1}(0)$ in a generic and close to the origin linear 
section of $f^{-1}(0)$,
inside a small ball.

\par In case (1), we consider the first non-trivial 
homotopy group 
for polynomials which may have certain singularities at infinity. This  
expands the class of polynomials which we can handle.
It appears to be useful when, in 
section \ref{nonisolated}, we prove the aforementioned result for  
 polynomials with one atypical value (say $0$):
the first non-trivial higher homology group of the generic fiber can 
be identified with the first non-trivial 
higher homotopy group of the complement
to the hypersurface $H \cap f^{-1}(0)$
in a generic linear subspace $H$ 
(the codimension of $H$ is the dimension of the singular locus of $f^{-1}(0)$).  
This identification immediately yields several conditions for vanishing 
and non-vanishing of the homology of generic fibers of such polynomials 
by applying vanishing and non-vanishing results for the 
homotopy groups of the complements (cf. section 
\ref{nonisolated}). These results are based on divisibility theorems 
for the orders of the homotopy groups (cf. \cite{Li})
and on vanishing results for hypersurfaces having mild singularities.

\par In the case of germs of holomorphic functions, we show that most of the results
discussed above in the case of polynomials with one atypical value can be 
proven in the local case. Besides the already mentioned 
relation between the homology of 
the Milnor fiber in certain 
dimension and the homotopy group of the complements to hypersurfaces in a ball,
we prove some divisibility results similar to those in the case
of polynomials with one atypical value. 
We suspect that this analog can be extended further and that 
the latter homotopy groups
can be related to geometry of singularities in a precise way, but because of 
technical difficulties, we postpone the discussion to later publication. 

\par The contents of this note are the following. In the next section we discuss
a way to measure the dimension of {\em singularities of a polynomial at infinity} (in a certain strong sense)
and a method of constructing such polynomials, motivated by \cite{Ti}.
In section \ref{higherhomotopy} we consider the homotopy groups of 
the complements
to hypersurfaces which may have singularities at infinity, 
expanding the results from \cite{Li}. 
We consider separately two situations. 
One is the case when the hypersurface is a generic 
fiber of a polynomial and another is the case of $f=t$, 
where $t$ is an atypical value.  
In the last section we 
prove results on the homology of smoothings and homotopy groups of 
the complements and the consequences discussed earlier in this
introduction in cases (1), cf. Theorem \ref{t:alex}, and (2), 
cf. Theorem \ref{t:local}.

%%%%%%%%%%%%%%%%%%%%%%%%%
%%%%%%%%%%%%%%%%%%%%%%%%%
%%%%%%%%%%%%%%%%

\section{Classes of hypersurfaces with singularities at infinity}

Let $V$ be a hypersurface in $\bC^{n+1}$. 
In \cite{Li}, the first author investigates the homotopy of the complement 
$\bC^{n+1}\setminus V$, which depends on the singularities of $V$ and also 
on {\em singularities at infinity of $V$}, defined as follows in {\em loc.cit.}:
\begin{equation}\label{eq:1} 
\Sing^\ity (V) := \Sing (\bar V \cap H^\ity),
\end{equation}
 where 
$H^\ity\subset \bP^{n+1}$ is the hyperplane at infinity and $\bar V \subset 
\bP^{n+1}$ is the projective closure of $V$.

The hypersurface $V$ is a fibre of a 
polynomial function $f: \bC^{n+1} \to \bC$. We show that the consideration 
of 
singularities of $f$ instead of those of $V$ may refine the study of the 
homotopy of the complement $\bC^{n+1} \setminus V$. Of course, we have to take 
into account the {\em singularities of $f$ at infinity}; our definition
extends Definitions 1.1, 2.3 in \cite{Ti} and has 
common flavor with Definition 2.2 in \cite{Li-duke}.
  
 %%%%%%%%%%%%%%   

  Consider an embedding of $\bC^{n+1}$ into some complex space $\bX$ such
 that 
there exists a proper algebraic morphism $f^{\bX} : \bX \to \bC$ extending $f$. In particular, the space $\bX$ compactifies the fibres of $f$. Let 
$\cS$ denote some Whitney stratification of $\bX$ such that $\bC^{n+1}$ is contained in a stratum.
  Let $\Sing f^{\bX}_{|\cS_i}$ denote the singular locus of the restriction of $f^{\bX}$ to $\cS_i$ (i.e. $\{ x\in \cS_i \mid \grad f^{\bX}_{|\cS_i}(x)=0 \}$). Then $\Sing_\cS f^{\bX} := \cup_{\cS_i\in \cS} \Sing f^{\bX}_{|\cS_i}$ is the {\em singular locus} of $f^{\bX}$ with respect to the 
stratification $\cS$. The Whitney conditions imply that $\Sing_\cS f^{\bX}$ is a closed analytic set. Notice 
that $\Sing_\cS f^{\bX}$ depends on the choice of the embedding $\bC^{n+1} 
\subset 
\bX$, whereas its intersection with $\bC^{n+1}$ does not,
 namely it is $\Sing f$.

\begin{definition}
\label{d:sing}
We call the germ of $\Sing_\cS f^{\bX}$ at the set $\bX^\ity :=\bX \setminus \bC^{n+1}$  the {\em 
singular set of $f$ at infinity} with respect to the proper
extension $f^{\bX}$ and to the stratification $\cS$. We say that {\em $f$ has isolated singularities at infinity} if the dimension of $\Sing_\cS f^{\bX}$ at any point of $\bX^\ity$ is $\le 0$. 
\end{definition} 

 %%%%%%%%%%%%%%
 Observe that, if $f$ has isolated singularities at infinity for some proper extension $f^{\bX}$ and some stratification $\cS$, then $\dim \Sing f \le 0$.
 
 We show that we have good control of the topology over a class of 
hypersurfaces which satisfy Definition \ref{d:sing}, although they might have a large singular locus at infinity in the 
sense of (\ref{eq:1}). 

%%%%%%%%%%%%%%

 Using the fact that $f^{\bX}$ is algebraic and Thom's First Isotopy Lemma, 
 it follows that the image $f^{\bX} (\Sing_\cS f^{\bX})$ is a finite set 
$\Lambda \in \bC$ (which depends on the embedding $\bC^{n+1} 
\subset \bX$) and that $f: \bC^{n+1}\setminus f^{-1}(\Lambda) \to 
\bC\setminus \Lambda$ is a locally trivial fibration (see e.g. \cite{Ve}). We assume that $\Lambda := \{ b_1, \ldots , b_r\}$ is minimal with this property and we call it the set of {\em atypical values}. We call $f^{-1}(b_i)$ an {\em atypical fiber}. It follows from the definition  that $\Lambda$
contains the critical values of $f$ and is contained in $f^{\bX} (\Sing_\cS f^{\bX})$.

 The embedding $\bC^{n+1} 
\subset \bX$ and extension $f^{\bX} : \bX \to \bC$ that we mainly use in this paper are the following.
%%%%%%%%%%%%%%%%%%%%%%%%%

\begin{example}\label{e:sing}
 To each coordinate $x_i$ of $\bC^{n+1}$  we associate a 
positive 
weight $w_i$ and write $f = f_d + f_{d-k} + 
\cdots$ where $f_j$ is the degree $j$ weighted-homogeneous part of $f$ and 
where $f_{d-k} \not= 0$.
 Let $\tilde f$ be the degree $d$ 
homogenization of $f$, in the weighted sense, by a new variable $z$ of weight 1, and define:
\begin{equation} \label{eq:x}
 \bX := \{ \tilde f(x,z) - tz^d =0\} \subset \bP(w) \times \bC,
\end{equation}  
where $\bP(w)$ denotes the weighted projective space $\bP(w_0, 
\ldots , w_n, 1)$. Since $\bP(w)$ is the space of orbits by the $\bC^*$-action on $\bC^{n+2}\setminus \{ 0\}$ given by $\lambda*x = (\lambda^{w_0}x_0, \cdots \lambda^{w_n}x_n, \lambda z)$, it has a canonical Whitney stratification by the orbit type (see \cite{Fe} or \cite[p. 21]{GLPW}), which is moreover the coarsest one. We take on $\bP(w)\times \bC$ the product stratification. This induces a stratification on the subspace $\bX 
\subset \bP(w)\times \bC$ and we consider the coarsest 
Whitney stratification $\cS$ on $\bX$ containing it. 
We then define $f^{\bX} : \bX \to \bC$ to be the 
projection on the second factor.  We identify $\bC^{n+1}$ with $\bX\setminus \bX^\ity$, where $\bX^\ity := \bX\cap \{ z=0\}$  denotes the ``divisor at 
infinity''of $\bX$.
Notice that the top-dimensional stratum of $\cS$ is $\bX \setminus \Sing \bX$ and it contains $\bC^{n+1}$, under the above identification.
\end{example} 
%%%%%%%%%%%%%% 
   
In Example \ref{e:sing}, if all the weights are 1, the singularities at infinity of $f$ can be estimated by the following 
practical criterion (the case $s=0$ has been proved in  \cite{Ti-lef} for any weights). We state the result only in case of weights 1 but this holds in fact for any weights, with a slightly more extended proof. 
%%%%%%%%%%%%%%%%%%%%%%%%% 

\begin{proposition}\label{p:sing}
Let $f : \bC^{n+1}\to \bC$ be a polynomial and consider the extension $f^{\bX} : \bX \to \bC$  as in Example \rm \ref{e:sing}. \it Let $\Sigma := \{ \grad f_d =0, f_{d-k}=0\} \subset \bP^{n+1}\cap 
\{ z=0\}$.
If the singular locus of $f$ in $\bC^{n+1}$ is of dimension $\le s$ and if $\dim \Sigma \le s$, then $\dim 
\Sing_\cS f^{\bX} \le s$.
\end{proposition}
%%%%%%%%
\begin{proof}
Since 
the singularities of $f^{\bX}$ on $\bC^{n+1} = \bX \setminus \bX^\ity$ are of dimension $\le s$ by 
hypothesis, we only have to look at singularities of $f^{\bX}$ within 
$\bX^\ity =\{f_d =0\} \times \bC \subset 
(\bP^{n+1} \cap \{ z=0\})\times \bC$. 
So we need to know the stratified structure of $\bX$ in the neighbourhood of $\bX^\ity$.
We prove in the following that $\bX^\ity \cap \Sing_\cS f^{\bX} \subset \Sigma \times \bC$.

Consider the hypersurface $\tilde \bX$ within $(\bC^{n+2}\setminus \{ 0\}) \times \bC$, defined by the same equation as in (\ref{eq:x}). 
 Notice first that the subset $\tilde \bX \cap \{ z=0\} \cap (\{ f_d =0 \} \setminus \{ \grad f_d = 0\})$ is contained in the regular part of  $\tilde \bX$. Next, at some point $\xi = (q, 0, t_0)\in \tilde \bX \cap \{ z=0\} \cap (\{ \grad f_d =0\}\setminus \{ f_{d-k} = 0\})$, we claim that there exists, locally at that point, a stratification of $\tilde \bX$ with the property that its strata are product-spaces by the $t$-coordinate (i.e., if $(p,t)$ belongs to a stratum, then $(p,t')$ belongs to the same stratum for all $t'\in \bC$).
 Indeed, the local equation of $\tilde \bX$ is $\tilde f(x,z) - tz^d = 0$ and can be written as
 $f_d(x) + z^k g =0$,
where $g(x,z,t) = f_{d-k}(x) + z h(x,z,t)$. Since  
$f_{d-k}(q) \not= 0$, one can define, locally at $\xi$, a new coordinate $z'= z\sqrt[k]{g}$, by choosing a $k^\th$ root
 of $g$. It follows that, locally,
 our hypersurface is equivalent, via an analytic change of 
coordinates at $\xi$, to the product of $\{f_d(x) + (z')^k = 0\}$ by the $t$-coordinate. Consequently, there exists  a local Whitney stratification at $\xi$ which is a product by the $t$-coordinate.
Notice that $\{z' =0\}$ corresponds to $\{z =0\}$ at $\xi$, hence the complement of $\{z' =0\}$ is nonsingular too. 

It follows that $\tilde \bX$ may be endowed with a global Whitney stratification $\cW$, such that $\tilde \bX^\ity$ is a union of strata and that, locally at each point $\xi \in \tilde \bX^\ity \setminus (\tilde\Sigma \times \bC)$, the strata contained in $\tilde \bX^\ity$ are products by the $t$-coordinate. In particular, the projection to $\bC$ is a stratified submersion at all points of $\tilde \bX^\ity \setminus (\tilde\Sigma \times \bC)$. Moreover, since in charts $\{ x_i =1\}$ the $\bC^*$ action on $\tilde \bX$ is the identity, it also follows that the projection $\bX \to \bC$, which is just our $f^\bX$, is a stratified submersion at every point $\xi \in \bX^\ity \setminus (\Sigma \times \bC)$.
 Taking the coarsest stratification on $\bX$  such that $\bX^\ity$ is a union of strata, the map  $f^\bX$ clearly remains a stratified submersion at all points of $\bX^\ity \setminus (\Sigma \times \bC)$. 

Our claim that $\bX^\ity \cap \Sing_\cS f^{\bX} \subset \Sigma \times \bC$,
is now completely proved. To conclude the proof of \ref{p:sing}, we notice that the stratified singularities of the restriction of $f^{\bX}$ to $\Sigma \times \bC$ can only occur on a finite number of fibres of $f^{\bX}$. Since
$\dim \Sigma \le s$, it 
follows that $\dim \bX^\ity \cap \Sing_\cS f^{\bX} \le s$. 
\end{proof}
%%%%%%%%%%%%%%%%%%%%%%%%%%%%%%%%%%%%%%%%%%%%%%%%%
%See the Examples in the next section.
\begin{corollary}\label{c:dim}
If $\dim \Sigma \le s$ then $\dim \Sing_\cS f^{\bX} \le s+1$ and, in particular, $\dim \Sing f \le s+1$.
\end{corollary}
\begin{proof}
The proof of \ref{p:sing} shows that, if $\dim \Sigma \le s$ then $\dim \bX^\ity \cap \Sing_\cS f^{\bX} \le s$. This implies that $\dim\Sing_\cS f^{\bX} \le s+1$. In turn, this yields $\dim \Sing f \le s+1$, since $\bC^{n+1}$ is a stratum of the stratification $\cS$ of $\bX$.
\end{proof}

%%%%%%%%%%%%%%%%%%%%%%%%%%%%%
%%%%%%%%%%%%%%%%%%%
%%%%%%%%%%%%%%%%%%%%%%%%%%%%
%%%%%%%%%%%%%%%%%%%
\section{Higher homotopy groups}
\label{higherhomotopy}

Let $V\subset \bC^{n+1}$ be a hypersurface. It may be a general fiber of a polynomial $f$, or an 
atypical one. In case $V$ is a general fiber, we have the following result on the homotopy type of the complement:
 %%%%%%%%%%%%%%

\begin{proposition}\label{prop:gen} 
 Let $f: \bC^{n+1} \to \bC$ be any polynomial and let $V$ be a general 
fibre of $f$. Then $\bC^{n+1}\setminus V \stackrel{\h}{\simeq} S^1 \vee 
S(V)$, where $S(V)$ denotes the suspension over $V$.

In particular, the cup-product in the cohomology ring of $\bC^{n+1}\setminus V$ is trivial. 
\end{proposition}
%%%%%%%%%%%%%%%%%%%%
\begin{proof}
Let $V = f^{-1}(\beta)$ and take a small enough closed disc $D\subset \bC$ 
centered at $\beta$. Take a path $\gamma_i$ from some point 
$\alpha$ on $\partial D$ to a small enough disc $D_i$ centered at an atypical value $b_i$ of $f$. Now, $f$ is a trivial fibration over $D$, hence $f$ is trivial over 
$\partial D$ too. Since $\bC^{n+1}\setminus V\stackrel{\h}{\simeq} 
f^{-1}(\cup_i (\gamma_i \cup D_i)) \cup f^{-1}(\partial D)$, it follows 
that 
$\bC^{n+1}
\setminus V$ is obtained from $f^{-1}(\partial D)\stackrel{\h}{\simeq} 
\partial D\times V$ by attaching the space  $f^{-1}(\cup_i (\gamma_i \cup 
D_i))\stackrel{\h}{\simeq}\bC^{n+1}$ over 
$f^{-1}(\beta)\stackrel{\h}{\simeq} V$. This is the attaching of a cone 
over $\{\beta\} \times V$ to $\partial D \times V$, so we get the claimed 
result.
\end{proof}

%%%%%%%%%%%%%%%%%%%%%%%%%
We may derive the following consequence; this has been proved in \cite[4.5]{Ti-lef} for a particular extension. \begin{corollary}\label{c:bou}
If $V$ is a general 
fibre of $f$ and $f$ has isolated singularities at infinity for some extension $f^{\bX}$, then $\bC^{n+1}\setminus V \stackrel{\h}{\simeq} S^1 \vee \bigvee_\lambda S^{n+1}$.
\end{corollary}
\begin{proof}
Since $f$ has isolated singularities at infinity, we may 
use a result of the second author \cite[Theorem 4.6, Corollary 4.7]{Ti} 
which works in our general setting. 
This bouquet theorem says that the general 
fibre of our $f$ is homotopy equivalent to a bouquet of spheres $\bigvee_\lambda S^n$, where $\lambda$ is the sum of the local Milnor numbers at singular points of $f$ and at singularities at infinity of $f^{\bX}$.
\end{proof}
%%%%%%%%%%%%%%
Corollary \ref{c:bou} holds in particular
if $f$ has isolated singularities in $\bC^{n+1}$ and $\dim \Sigma \le 0$,
 by Proposition \ref{p:sing}. 

%%%%%%%%%%%%%%
\begin{example}\label{ex:1}
 $f: \bC^3 \to \bC$, $f = x+ x^2y + z^2$. \\
 This polynomial has no singularities in 
$\bC^3$. 
If we consider the extension $f^{\bX}$ as in Example \ref{e:sing}, with all weights equal to 1, then $f$ has isolated singularities at infinity since it satisfies the assumptions of Proposition \ref{p:sing} for $s=0$. According to Proposition \ref{prop:gen}, 
we have the homotopy equivalence $\bC^3 \setminus V \stackrel{\h}{\simeq} 
S^1 \vee S^3$; in particular $\pi_2(\bC^3 \setminus V)=0$.
    This works for general fibres, i.e. for $V=f^{-1}(t)$, $\forall t\not= 0$, since the only atypical value  of $f$ is $0$, as one can easily check.
\end{example}

%%%%%%%%%%%%%%

\begin{example}\label{ex:2}
 $f: \bC^4 \to \bC$, $f = x_1^4x_2^4 + (x_1 + x_2)^6 + x_3^5 + x_4^4 + 
x_1^2$.\\
 Note that, according to definition 
(\ref{eq:1}), $V$ is not transversal at infinity along a 2-dimensional set. 
Nevertheless, we may observe that $f= g(x_1, x_2) + h(x_3, x_4)$ is a sum of two polynomials in separate variables. For $g$, we have that 
$\dim \Sing g \le 0$ and that there are 
no 
singularities at infinity, since $\Sigma = \emptyset$ (use Proposition \ref{p:sing}). We get that the general fibre of $g$ is, homotopically, a bouquet  $\bigvee S^1$.
On the other hand, the polynomial $h$ is weighted homogeneous with a unique singularity at the origin, hence its general fibre is $\stackrel{\h}{\simeq}\bigvee S^3$.
 
  By a Thom-Sebastiani result, the general fibre of $f= g+h$ has the homotopy type of a bouquet $\bigvee S^3$. 
  Now, by Proposition \ref{prop:gen}, for a general fibre $V$ of $f$, 
the 
complement $\bC^4 \setminus V$ is homotopy equivalent to $S^1 \vee \bigvee
S^4$.  
\end{example}

%%%%%%%%%%%%%%
For the complement of a hypersurface $V$ which is an atypical fibre of a polynomial $f$, we have the following result:
 %%%%%%%%%%%%%%

\begin{proposition}\label{p:spec} 
Let $V = f^{-1}(0)$.
 If the general fibre of the polynomial function $f: \bC^{n+1} \to \bC$ is 
$s$-connected, $s\ge 2$, then $\pi_i(\bC^{n+1}\setminus V) =0$, for $1<i\le 
s$, and $\pi_1(\bC^{n+1}\setminus V) =\bZ$.
   In particular, if $f$ has isolated singularities at infinity in some extension $f^{\bX}$, then $\pi_i(\bC^{n+1}\setminus V) =0$ 
  for $1<i\le n-1$.
\end{proposition}
\begin{proof} 
  We use the notations $D_i$, $\gamma_i$ as in the proof of Proposition \ref{prop:gen} and take $b_1 := 0$, so that $V = f^{-1}(b_1)$. We first claim that the space $T_1 := 
f^{-1}(\cup_{i\not= 1} \gamma_i \cup D_i)$ is homotopy equivalent to a general fibre 
$F := f^{-1}(\alpha)$ to which one attaches cells of dimension $\ge s+2$, in other words that the pair $(T_1,F)$ is $(s+1)$-connected. 

Since $F$ is 
$s$-connected by hypothesis, we have that $(\bC^{n+1}, F)$ is $(s+1)$-connected.
Then, by excision in homology, we have that  $H_j(\bC^{n+1}, F) = H_j(T_1,F) \oplus H_j(f^{-1}(D_1),F)$, for any $j$. Hence $H_j(T_1,F) =0$ for $j\le s+1$.
By Blakers-Massey theorem \cite{BM},  the excision works in homotopy within a certain range. Namely, since $F$ is $s$-connected, we get that the inclusion:
\[ \pi_j(T_1,F) \oplus \pi_j(f^{-1}(D_1),F) \to \pi_j(\bC^{n+1}, F) \]
is an isomorphism for $j\le s-1$ and an epimorphism for $j=s$.
This shows that $(T_1,F)$ is $(s-1)$-connected and in particular simply connected, since $s\ge 2$.
We may furthermore apply the relative Hurewicz isomorphism theorem and get that  $\pi_j(T_1,F)$ is trivial for $j\le s+1$. (We also get the isomorphism $\pi_{s+2}(T_1,F) \simeq H_{s+2}(T_1,F)$.)
  By Switzer's result \cite[Proposition 6.13]{Sw}, it follows that $T_1$ is 
homotopy equivalent to the space $F$ to which one attaches cells of 
dimension $\ge s+2$. The claimed property is proved. 
 
 Next, we have that $\pi_j (f^{-1}(\partial D_1)) =0$ for $1< j\le s$ and that $\pi_1(f^{-1} (\partial D_1) =\bZ$, due to the homotopy exact sequence of the fibration $f_| : f^{-1}(\partial D_1) \to 
\partial D_1$ and the $s$-connectivity of the fibre. (Note that this holds even for $s=1$.)

  Finally, $\bC^{n+1} \setminus V$ is obtained from $f^{-1}(\partial D_1)$ by attaching the space $T_1$ over a general fibre $F$, which, we have proved above, means attaching only cells of dimension $\ge s+2$. It follows that $\pi_j 
(\bC^{n+1}\setminus V) =0$ for $1< j\le s$ and that $\pi_1(\bC^{n+1} \setminus V) =\bZ$.
  
  The second statement follows from the first one. Indeed, if $f$ has isolated singularities at infinity then,  as mentioned in the proof of Corollary \ref{c:bou}, the general fibre of $f$ is $(n-1)$-connected.
\end{proof}

%%%%%%%%%%%%%% 

\begin{theorem}\label{p:connec}
Let $f: \bC^{n+1} \to \bC$ be a polynomial function and let's consider the embedding $\bC^{n+1} \subset \bP(w_0, \ldots , w_n, 1)$, for some system of weights $w$, as in
Example \ref{e:sing}. Suppose that $\dim \Sing_\cS f^{\bX} \le k$. Then:
\begin{enumerate}
\rm \item \it for a general fibre $F$
of $f$, $\pi_i (\bC^{n+1}\setminus F) = 0$, for $2\le i \le n-k$.
\rm \item \it  for an atypical fibre $V$
of $f$, $\pi_i (\bC^{n+1}\setminus V) = 0$, for $2\le i\le n-k-1$.
\end{enumerate}
\end{theorem}
\begin{proof}
If we prove that the general fibre $F$ is $(n-k-1)$-connected, then (a) 
follows by 
Proposition \ref{prop:gen} and (b) follows by  Proposition \ref{p:spec}.

  So let us show that the condition  $\Sing_\cS f^{\bX} \le k$ indeed implies that the 
general fibre $F$ is $(n-k-1)$-connected.
This is true in the particular case when all the weights are $1$,
by \cite[Theorem 5.5]{Ti}. The proof in {\em 
loc.cit.} goes by induction and uses generic hyperplane sections in $\bP^{n+1}$, which do not exist in the case of weighted projective space.

Nevertheless, the proof could work in a similar spirit,
 provided that we can use, instead of generic hyperplane sections, a class
of hypersurface sections with good enough properties. Let us start defining that.

 Consider the finite map $\Psi : \bC^{n+1} \to \bC^{n+1}$, given by $(x_0, \ldots 
, x_n) 
\mapsto   (x_0^{m_0}, \ldots , x_n^{m_n})$, where $m_i =N/w_i$ and $N$ is a 
common 
multiple of all $w_i$, $i=\overline{0,n}$.  This induces $\bar \Psi : 
\bP(w_0, \ldots , w_n) \to \bP^n$ and so a finite map $\bX \subset \bP(w)\times \bC \stackrel{\hat \Psi}{\to} \bP^{n+1}\times \bC$. The class of ``generic'' hypersurfaces
will be an open subset of $\{ a_0x_0^{m_0} +\cdots + a_nx_n^{m_n} = s \mid a_i, s \in \bC \}$. Actually, we prove our statement by reduction to the space $\bP^{n+1}$, via the finite map $\hat \Psi$.

 We denote by  $H_s$ the affine hyperplane $\{ l_H =s\} \subset \bC^{n+1}$, 
where
$H\in \check\bP^{n}$ is a hyperplane defined by a linear form $l_H : 
\bC^{n+1} \to 
\bC$. We consider the restriction of $f$ to $\cH_s := \Psi^{-1}(H_s)$.
We denote by $\cS_{\cH_s}$ the coarsest  Whitney stratification 
of the space $\bH_s := \bX \cap (\bar \cH_s \times \bC)$,
  where  $\bar \cH_s$ denotes the degree $N$ weighted projective 
hypersurface 
$\{ l_H\circ \Psi (x) - sz^N =0\}$.

By eventually refining the stratification $\cS$, we may assume without loss of generality that the restriction of $\hat\Psi$ to each stratum is an unramified covering. To the images by $\hat\Psi$ of the strata and of the levels of $f^{\bX}$, we may apply the method of slicing by generic hyperplanes, as described in \cite[\S 5]{Ti}. Then we may transfer back the transversality results via $\hat\Psi^{-1}$.

In this way, by using \cite[Lemma 5.4]{Ti}, it follows that there exists a Zariski-open set 
$\Omega\subset \check \bP^{n}$
 and a finite set $A\subset \bC$ such that, if
$H\in \Omega$ and $s\in \bC \setminus A$ and if $\dim \Sing_\cS f^{\bX} \ge 
1$, then
$\dim \Sing_{\cS_{\cH_s}}(f_{| \cH_s})^{\bH_s} \le \dim \Sing_\cS f^{\bX} -1$.

  Slicing a general fibre $F$ of $f$ by a ``generic" 
hypersurface $\cH_s$ gives a general fibre of the restriction $f_{|\cH_s}$. 
Moreover, 
by the Lefschetz type theorem, the pair $(F, F\cap 
\cH_s)$ is 
$\dim F -1$ connected.
  
 The general setting of \cite{Ti} and the use of the pull-back by $\hat\Psi$ allow one to continue this slicing procedure
 until the singularities of the restriction of $f$ become zero-dimensional.
 Namely, there exist $k$ generic hyperplanes $H^i \in \check\bP^n$ and 
generic 
 $s_1, \ldots , s_k \in \bC$
 such that $A := \cH^1_{s_1}\cap \cdots \cap \cH^k_{s_k}$ is the global 
Milnor 
fibre of a weighted homogeneous affine complete intersection with isolated 
singularity at the origin
 and that the restriction $f_| :  A \to \bC$ has isolated singularities in 
the 
affine and at infinity, in the sense of Definition \ref{d:sing}.
 
  In this situation, we may apply to $f_| :  A \to \bC$ the results 
concerning 
isolated singularities at infinity \cite[Theorem 4.6 and Corollary 
4.7]{Ti}, namely: $A$ is obtained from a general fibre of $f_{|A}$ by attaching a finite number of cells of dimension $n+1-k$.
   Since $A$ is homotopy equivalent to a bouquet of spheres of dimension 
$n+1-k$, it follows that this general fibre is homotopy 
equivalent to a bouquet of $n-k$ spheres, since it is a Stein space of dimension $n-k$ and $(n-k-1)$-connected. By tracing back the vanishing of homotopy in the  slicing 
sequence,
  we get that the general fibre $F$ is at least $(n-k-1)$-connected.  
\end{proof}

%%%%%%%%%%%%%%%%%%%%%%%%%%% 
 The assumption of Theorem \ref{p:connec} holds in particular if $\dim (\Sigma \cup \Sing f) \le k$, by Proposition  \ref{p:sing}. About $\pi_1$, we can say the following:
 %%%%%%%%%%%%
\begin{remark}\label{r:connec} 
\begin{enumerate}
\item In case of a generic fibre $V$, it follows from the proof of Theorem \ref{p:connec}(a) via Proposition \ref{prop:gen}, that $\pi_1(\bC^{n+1} \setminus V)$ is trivial as long as $n-k\ge 1$, where $k\ge 0$.

\item In case of an atypical fibre $V$, one may prove that $\pi_1(\bC^{n+1} \setminus V)$ is trivial, as long as 
$n-k\ge 2$, where $k\ge 0$. We have to modify the arguments in the proofs of the above results, as follows. We slice as in the proof 
of Theorem \ref{p:connec}, but instead of following a generic fibre $F$, we work with the atypical $V$. We get that the restriction $f_| :  A \to \bC$ has isolated singularities in the  affine and at infinity
and that $V\cap A$ is an atypical fibre of $f_{|A}$.
Next, revisiting the proof of Proposition \ref{p:spec}, we notice that 
one can prove by some different arguments that, for our restriction $f_{|A}$, the pair $(T_1, F)$ is $(n-k)$-connected. Namely, we may use here again \cite[Theorem 4.6 and Corollary 4.7]{Ti}, as in the proof of Theorem \ref{p:connec}, to show that $T_1$ is obtained from a general fibre of $f_{|A}$ by attaching a number of cells of dimension $n+1-k$.
 The last part of the proof of Proposition \ref{p:spec} still works in case $s=1$. 
\end{enumerate}
\end{remark}

\begin{note}\label{n:dp}
 The proof of Theorem \ref{p:connec} yields, in particular, that the general fibre of $f$ is $(n-k-1)$-connected and that any atypical fibre of $f$ is at least $(n-k)$-connected. This bound for the connectivity of the fibres appears to be sharp.
\end{note} 
%%%%%%%%%%%%%%%%%%%%%%%%%%%%%%%
%%%%%%%%%%%%%%%%%%%%%%%%%%%%%%%
%%%%%%%%%%%%%%%%%%%%%%%%%%%%%%%
\section{Monodromy of non-isolated singularities} 
\label{nonisolated}
%%%%%%%%%%%%%% 

\subsection{The global case}

We show how the monodromy of certain non-isolated 
singularities is related to the ``first'' non-vanishing homotopy group of the complement of the  hypersurface singularity. 
Though most of the material below can be carried out in the framework 
of definition \ref{d:sing} we shall start by working with the weighted projective embedding 
 $\bC^{n+1} \subset \bX \subset \bP^{n+1}(w)\times \bC$, as defined at \ref{e:sing}. Let us assume that $\dim \Sing_\cS f^{\bX} 
\le k$ and denote by $V= f^{-1}(0)$ an atypical fibre. The generic fiber $F$ of $f$ is $(n-k-1)$-connected, by Theorem \ref{p:connec}. Our goal is to calculate
the monodromy acting on {\it the first possibly non-trivial group} $H_{n-k}(F)$. 
\par Let $L_{k-1} = \cH_1 \cap \cdots \cap \cH_{k-1}$ denote the intersection of $k-1$ generic hypersurfaces, as defined in the proof of Theorem \ref{p:connec}.
Then the Lefschetz theorem yields:
\begin{equation} \label{lefschetzcut}
H_{n-k}(F)=H_{n-k}(F \cap L_{k-1}),
\end{equation} 
where $F \cap L_{k-1}$ can be 
viewed as a generic fiber of the polynomial $f_{|L_{k-1}}$. Moreover, $\dim \Sing_{\cS'} (f_{|L_{k-1}})^{\bL_{k+1}} \le 1$, after \cite{Ti}, where  $\bL_{k+1} := \bX \cap (\bar L_{k-1}\times \bC)$ (see the proof of \ref{p:connec} for the definition of the stratification $\cS'$ on $L_{k-1}$).
 
\par On the other hand if $L_k = L_{k-1}\cap \cH_k$ is the cut by another generic hypersurface, then 
by the Zariski-Lefschetz theorem of Hamm and L\^e \cite{HL} we have:
\begin{equation}\label{lefschetzhomotopy}
\pi_{j}(L_k\setminus L_k\cap V) \simeq \pi_{j}(\bC^{n+1}\setminus V),
\end{equation}
for $j\le n-k$. This is an isomorphism 
of $\bC[\bZ]$-modules. Moreover, for $i <n-k$, we have $\pi_1(L_k\setminus L_k\cap V)=\bZ$
and $\pi_i(L_k\setminus L_k\cap V)=0$, by \ref{p:spec}. The polynomial $f_{|L_k}$ has only isolated singularities
at infinity and $\pi_{n-k}(L_k \setminus L_k \cap V)$ is the {\it first possibly non-trivial 
homotopy group} (cf. Theorem \ref{p:connec}).
 
%%%%%%%%%%%%%%%%%%%%%%%%%%
\begin{definition}\label{d:delta}
Let 
\begin{equation}\label{eq:delta}
\pi_{n-k}(\bC^{n+1}\setminus V)\otimes \bQ = \bigoplus_i \bQ[t,t^{-1}]/(\lambda_i) \oplus \bQ [t,t^{-1}]^\kappa
\end{equation} 
be the decomposition as $\bQ [t,t^{-1}]$ module, where $\lambda_i$ are some polynomials, defined up to units in $\bQ [t,t^{-1}]$. 

When $\kappa = 0$ in the decomposition (\ref{eq:delta}), one calls the product $\prod_i \lambda_i$ the {\em order of} $\pi_{n-k}(\bC^{n+1}\setminus V)\otimes \bQ$. We denote it by 
$\Delta (\bC^{n+1}\setminus V)$ or by $\Delta(L_k \setminus V\cap L_k)$.
In case $n-k=1$, this is nothing else than the Alexander polynomial of the curve $L_k \cap V$ in $L_k$.
When $\kappa \not= 0$, one says the order is $0$. 
\end{definition}
%%%%%%%%%%%%%%%%%%%%%%%%%%

\begin{theorem}\label{t:alex}
Let $V= f^{-1}(0)$ be the only atypical fiber\footnote{For a partial classification of  polynomials in 2 variables with one atypical value, we refer
 to \cite[4.4]{Ti-reg}.} of 
a polynomial $f$ with $\Sing_{\cS} f^{\bX} \le k$. Let $L_k$ denote the intersection of $k$ generic hyperplanes in $\bC^{n+1}$.
Then
\[ \pi_{n-k}(L_k \setminus V \cap L_k)=H_{n-k}(F,\bC) \]
as $ \bC[\bZ]$-modules. In particular, for the orders of these modules we have:
\[ \Delta(L_k \setminus V\cap L_k) = \det[h_0 - t\cdot\id : 
H_{n-k}(F,\bC) \to 
H_{n-k}(F,\bC)],\]
where $h_0$ is the monodromy of the general fibre $F$ of $f$ around 
the value $0$.
\end{theorem}
%%%%%%%%%%%%%%%%%%%%%%%%

\begin{proof}
Since $f$ has at most one atypical value, the infinite cyclic cover
  $\widetilde{\bC^{n+1} \setminus f^{-1}(0)}$ is homotopy equivalent to 
the general fiber $F$.  The monodromy $h_0$ is the deck transform of the infinite cyclic cover
$\widetilde{\bC^{n+1} \setminus f^{-1}(0)} \rightarrow 
{\bC^{n+1} \setminus f^{-1}(0)}$. This, together with (\ref{lefschetzhomotopy}) 
yields the first equality above.
Moreover, we have the following:
 \[\begin{array}{ll} \Delta(L_k\setminus L_k\cap V) & =  \mbox{order of } 
\pi_{n-k}(L_k\setminus L_k\cap V) \otimes \bC = \\
 \  & = \mbox{order of } \pi_{n-k}(\bC^{n+1}\setminus V) \otimes \bC= 
\\
\  & = \mbox{characteristic polynomial of the deck transform on } H_{n-k} 
(\widetilde{\bC^{n+1}\setminus V}) =\\
\  & = \mbox{characteristic polynomial of the monodromy on } H_{n-k} 
(F, \bC).\end{array} \]
\end{proof}
%%%%%%%%%%%%%%%%%%%%%%%%%%%%%%%
The above theorem can be used to obtain results for the homology of Milnor fibers of polynomials with non-isolated singularities, as follows.

%%%%%%%%%%%%%%%% 4.3  %%%%%
\begin{corollary} \label{c:cor}
Under the assumptions of Theorem \ref{t:alex},
let $k=1$ and $n>2$. Then the number of cyclic factors  corresponding to 
$t-1$ in the cyclic decomposition of $H_{n-1}(F,\bQ)$ is 
equal to the dimension of $H_{n-1}(\bC^{n+1} \setminus V, \bQ)$.
In particular $\rk H_{n-1}(F,\bQ) \ge \rk H_{n-1}(\bC^{n+1}\setminus V,\bQ)$ 
and this
equality takes place if and only if $\Delta(\bC^{n+1}\setminus V)$ has 
no other roots except 1 and $k_i=1$ for any $i$.

If $n=2$ and $V$ is irreducible, then $\rk H_{n-1}(F,\bQ)=0$ if and only if $\Delta(\bC^{n+1}\setminus V) =1$.
\end{corollary}
\begin{proof} Let 
\[ H_{n-1}(F,\bQ)= H_{n-1}(\widetilde{\bC^{n+1}\setminus V},\bQ)=
\bigoplus_{i=1}^{l'} \bQ [t,t^{-1}]/(t-1)^{k_i} \bigoplus \oplus_{j=1}^{l''} \bQ [t,t^{-1}] /(\lambda_j), \]
where $\lambda_j$ does not have 1 as root.

From the spectral sequence $H_p(\bZ,H_q (\widetilde {\bC^{n+1}\setminus V}), \bQ)
\Rightarrow  H_{p+q}(\bC^{n+1}\setminus V, \bQ)$ and since
$H_{n-2} (\bC^{n+1} \setminus  V, \bQ) = 0$, we derive that 
$H_{n-1}(\widetilde {\bC^{n+1}\setminus V}, \bQ)^{Inv}=
H_{n-1}(\bC^{n+1}\setminus V, \bQ)$.
Since 
$H_{n-1}(\widetilde {\bC^{n+1}\setminus V}, \bQ)^{Inv}$ is the 
kernel of the multiplication by $t-1$ on 
$H_{n-1}(\widetilde {\bC^{n+1}\setminus V}, \bQ)$, the result follows
for $n>2$. If $n=2$, then the same spectral sequence shows that 
the number of cyclic summands is equal to $\rk H_1(\bC^3 \setminus V,\bZ) -1$ 
and the rank of the latter homology group can be identified with the number of 
irreducible components of $V$. 

\end{proof}

%%%%%%%%%%%%%%%%%%%%%%
%%%%%%%%%%%%%%%%%%%%%%%%%%%%%%%

\begin{example}
Let us consider a polynomial of the form 
\begin{equation}\label{eq:curve}
f=P_a(x,y,z)^b+ P_b(x,y,z)^a ,
 \end{equation}
 where $P_b$ and $P_a$ 
are generic and homogeneous of degrees $b$, resp. $a$. The homogeneity of $f$ 
yields that $0$ is the single atypical value; the singularities of $f$ form the 
union of lines corresponding to the points $\{ P_a=P_b=0\} \subset \bP^2$ and $\Sing_\cS f^{\bX}$ is just the closure of $\Sing f$ in $\bX$. The standard identification of 
the Milnor fiber of $f$ at $0$, denoted $M_f$, with the the cover of degree $ab$ of the complement of the projective curve given by $f=0$ 
and the calculation of the 
Alexander polynomial for such curves yields that 
$H_1(M_f,\bZ)=\bZ^{(a-1)(b-1)}$ (cf. \cite{Li-pr}). 
On the other hand, by \ref{t:alex}, the characteristic 
polynomial of this singularity is the Alexander polynomial 
of the affine curve $P_a(x,y,1)+P_b(x,y,1)=0$ (since the plane $z=1$ is 
generic), which is ${{(t^{a+b}-1)(t-1)} \over {(t^a-1)(t^b-1)}}$.   

More generally, we can consider the 
polynomial:
\begin{equation}
f=\sum P_{a_1 \cdots \hat a_i  \cdots a_n}(x_1,...,x_{n+1})^{a_i}
\label{hypersurface}
\end{equation} 
where $P_k$ denotes a generic homogeneous polynomial of degree $k$.
Again it satisfies the conditions of Theorem \ref{t:alex} with $k=1$. We have $H_{n-1}(M_f,\bC) 
\ne 
0$, according to \cite{Li-2}.
\end{example}

%%%%%%%%%%%%%%%%
%%%%%%%%%%%%%%%%%%%%%%%%%%%%%%%

\begin{example} One can obtain new examples of polynomials 
by applying automorphisms 
of $\bC^{n+1}$. For example, let  $f_{i,j}(x_j,..,x_{n+1}), 1 \le i <j \le n+1$ be arbitrary 
polynomials. 
The automorphism 
$$x_i \rightarrow x_i+ \sum_{j>i} f_{i,j}(x_j,...,x_{n+1})x_j$$
applied to the polynomial (\ref{hypersurface})
yields non-homogeneous examples of polynomials with non-isolated singularities
and non-trivial monodromy on the first non-vanishing homology of the Milnor fiber.  
\end{example}

%%%%%%%%%%%%%%%%
%%%%%%%%%%%%%%%%%%%%%%%%%%%%%%%%%%%%%%%%%%%%%%%%
\smallskip

 In the remainder we refer to the embedding of $\bC^{n+1}$ into $\bX\subset \bP^{n+1}\times \bC$, where $\bP^{n+1}$ is the usual projective space. 
Recall that $\Sing V \subset \Sing f \subset\Sing_\cS f^{\bX}$, where $V = f^{-1}(0)$.
  We consider the restriction $f_{|H}$ of the polynomial $f$ to a generic hyperplane in $\bC^{n+1}$. Let $\cS'$ denote the Whitney stratification induced by $\cS$ and the cut by $H$, as defined in \cite[5.2]{Ti} and recalled in \ref{p:connec} above. Therefore $\Sing_{\cS'} (f_{|H})^{\bH}$ is well defined.
  According to \cite[\S 5]{Ti}, there exist, and we shall use in the following, generic hyperplanes $H$ 
 such that $\bar V \cap \Sing_{\cS'}(f_{|H})^{\bH} = \bar H\cap\bar V \cap \Sing_\cS f^{\bX}$.
 
 Suppose  $\dim \Sing V =1$. For each irreducible component 
$\Sigma_i$ of $\Sing V$ let $\Delta_i$ be the characteristic 
polynomial of the monodromy of the isolated singularity one obtains as the transversal 
intersection of a small disk within $H$ with the one-dimensional stratum  $\Sigma_i$. We may call it {\it the horizontal} monodromy corresponding to $\Sigma_i$, in analogy to the case of germs, cf. \cite{St}.

 Suppose also that $\dim (\bX^\ity \cap \bar V \cap\Sing_\cS f^{\bX})  =1$ and denote by $\Sigma^\ity_i$
some one-dimensional irreducible component of $\bX^\ity \cap \bar V \cap\Sing_\cS f^{\bX}$. Then $\Delta_j$ denotes the characteristic 
polynomial of the monodromy of the isolated singularity at infinity
of the polynomial $f_{|H}$, for generic $H$.
Remark that the monodromy around an isolated singularity at infinity is well defined. Indeed, if $g : \bC^n \to \bC$ is a polynomial and if $p_j\in \bX^\ity$ is an isolated singularity of $g^\bX$, there is a locally trivial fibration
\[ B_j \cap \bX_{D^*}\setminus \bX^\ity \to D^*,\]
where $D$ and $B_j$ are Milnor data, i.e. $B_j$ is a small Milnor ball centered at $p_j$ and $D$ is a small enough disk centered at $g^\bX(p_j)$.

Following \cite[4.5-4.7]{Li}, let us denote by $\Delta^\ity$ the order of $\pi_{n-1}(S^\ity\cap H \setminus (S^\ity \cap V \cap H))$, for a sphere $S^\ity \subset \bC^{n+1}$ of sufficiently large radius.

With these preliminaries, we can prove the following.

%%%%%%%%%%%%%%%%%%% 4.6 %%%%%%%%%%%%

\begin{corollary}\label{c:div} 
Let $f: \bC^{n+1} \rightarrow \bC $ be a polynomial
such that $V= f^{-1}(0)$ is the only atypical fibre and such that
 $\dim \Sing_\cS f^{\bX} =1$ and let $F$ denotes a  generic fiber of $f$.
If none of the roots of $\Delta_i$'s and of $\Delta_j$'s 
distinct from $1$ is a root of $\Delta^\ity$ and $k_i=1$ for the cyclic 
summands in $H_{n-1}(F, \bQ)$ corresponding to $1$ then 
$H_{n-1}(F,\bQ)= H_{n-1}(\bC^{n+1}\setminus V, \bQ)$.
\end{corollary}

\begin{proof}
 Note that $f_{|H}$ does not have a single atypical value, in general, even if $f$ has only one. Due to the genericity of $H$,
the polynomial $f_{|H}$ has isolated singularities at infinity
and for such a polynomial the divisibility theorem \cite[4.3]{Li} works. More precisely, the cited result can be extended in case of isolated singularities at infinity in the sense of this paper. The conclusion from {\em loc.cit.} is that $\Delta(\bC^{n+1} \setminus V) = \Delta(H \setminus V\cap H)$ divides the the product of $\prod_i\Delta_i \cdot (1-t)^\kappa$, for some non-negative integer $\kappa$.

To complete the proof we have to use 
also \cite[Theorem 4.5]{Li}. This result applies to $f_{|H}$, which has isolated singularities at infinity. Actually, it can be extended to our new definition of singularities at infinity and it yields that $\Delta(H \setminus V\cap H)$ divides $\Delta^\ity$ of $f_{|H}$. Now our claim follows from Corollary 
\ref{c:cor}.
 \end{proof}

%%%%%%%%%%%%%%%%
Still in case $\dim \Sing V = 1$, our statement \ref{c:div} may become more precise, provided that we have a closer control on singularities at infinity.

\begin{corollary}\label{c:special}
 Let $f: \bC^{n+1} \rightarrow \bC $ be a polynomial
with one atypical fibre $V= f^{-1}(0)$ and such that
$\dim \Sing f  =1$ and that
$\dim \Sing (\bar V\cap H^\ity) =0$ (where $H^\ity$ denotes the hyperplane at infinity in $\bP^{n+1}$).
 
If none of the roots of $\Delta_i$'s and $\Delta_j$'s distinct from $1$ 
is a root of unity of degree $d = \deg f$, and $k_i=1$ for the cyclic 
summands corresponding to $1$ 
then  $H_{n-1}(F,\bQ)=H_{n-1}(\bC^{n+1}\setminus V, \bQ)$.
\end{corollary}
\begin{proof}
Due to the condition on singularities, the cut by a generic $H$ gives that 
$V\cap H$ is transversal at infinity $H^\ity$. We may apply \cite[Theorem 4.8]{Li}, which amounts to saying that $\Delta^\ity$ has only roots of unity of degree $d$. The conclusion then follows by 
 Corollary \ref{c:div}. 
\end{proof}
\begin{remark}\label{r:singinf}
Corollary \ref{c:special} remains true when replacing the condition $\dim \Sing (\bar V\cap H^\ity) =0$ by the more general one $\dim \bX^\ity \cap \Sing_\cS f^{\bX} =0$. See the preliminaries before Corollary \ref{c:div}.

\end{remark}
%%%%%%%%%%%%%%%%

%%%%%%%%%%%%%%%%%
\begin{example} Let $f: \bC^3 \to \bC$ be as in Corollary \ref{c:special} (e.g. homogeneous or obtained via an automorphism applied to a homogeneous 
polynomial) and having transversal $A_1$ or $A_2$-singularities 
along each of the strata of its singular locus. If the degree of $f$ 
is not divisible by $6$, then the rank of $H_1(F,\bC)$ is equal to the number
of irreducible components of the atypical fiber minus one. 
Indeed, the roots of the characteristic polynomial of the monodromy 
for an $A_2$ singularity are roots 
of unity of degree $6$. Hence if $d$ is not divisible by $6$ then the 
only root of the characteristic polynomial is $1$. 
Therefore, by \ref{c:cor}, the rank of $H_1(F,\bZ)$ is the equal to 
the rank of $H_1(\bC^3 \setminus f^{-1}(0))$ minus one 
and the rank of the last homology group  
is equal to the number of irreducible components of $f^{-1}(0)$.  
\end{example} 

Other vanishing, resp. non-vanishing results (cf. \cite{Li}) 
combined with \ref{t:alex} yield 
corresponding vanishing, resp. non-vanishing
results for the homology of the Milnor fibers of
non-isolated singularities.
Let us quote one of the results along these lines. 

%%%%%%%%%%%%%% 4.8

\begin{corollary} Let $f: \bC^3 \rightarrow \bC $ be a polynomial
having $0$ as single atypical value, such that $f^{-1}(0)$ is irreducible, that $\dim \Sing f^{-1}(0)=1$ and that $\dim (\Sing_\cS  f^\bX \cap \bX^\ity) =0$. 
Let the transversal type of singularities along the strata of $\Sing f^{-1}(0)$
be either $A_1$ or $A_2$. Assume that the degree $d$ of $f$ is divisible by $6$.
Let $\Sigma = \Sing f^{-1}(0) \cap H$, for some generic hyperplane $H$ and let the superabundance of the 
curves in $H$ having the degree equal to $d-3-{d \over 6}$ and passing through $\Sigma$ be $s$.
Then the characteristic polynomial of the monodromy of $f$ acting on $H_1(F,\bC)$ 
is equal to $(t^2-t+1)^s$. 
\fin
\end{corollary}
%%%%%%%%%%%%%%%%%%%%%%%%%%%

\subsection{The local case}

The above global case, with a single atypical value, is somehow close to a semi-local situation. So, let us formulate the problem 
in the local case.

Let $g : (\bC^{n+1}, 0) \to (\bC, 0)$ denote a germ of a holomorphic 
function, with 
$\dim \Sing g \le 1$. Then $g$ has a $(n-2)$-connected Milnor fibre $M_g$, by  
Kato and Matsumoto's result \cite{KM}. 
Take a closed Milnor ball $B$ at $0$, with boundary $S := \partial B$ and 
take a general hyperplane $H$ passing close to the origin. Denote $B_H = 
B\cap H$, $S_H = S\cap H$ and $V = g^{-1}(0)$. Then $B_H \cap V$ has at most isolated singularities, 
since $\dim \Sing g \le 1$. Let $\{ \Delta_i\}_i$ be the collection of characteristic polynomials of the monodromies of the singularities of $B_H \cap V$.  With these notations, we have the following:
%%%%%%%%%%%%%%%%
%%%%%%%%%%%%%%%%%%
%%%%%%%%%%%%%%%%

\begin{theorem} 
\begin{enumerate}
\rm \item \it $\pi_i (S_H \setminus S_H\cap V) =
\pi_{i}(B_H \setminus B_H \cap V)
= 0$, for $1< i \le n-2$ and 
$\pi_1(S_H \setminus S_H\cap V) =\pi_{1}(B_H \setminus B_H \cap V)=
\bZ$.
\rm \item \it $\pi_{n-1} (S_H \setminus S_H\cap V)$ is a torsion 
$\pi_1$-module.
\rm \item \it The module 
$\pi_{n-1}(B_H \setminus B_H \cap V) \otimes \bC$ is
a $\bC[\bZ]$-torsion module 
and the order of $\pi_{n-1}(B_H \setminus B_H \cap V) \otimes \bC$ 
divides the 
order 
of $\pi_{n-1} (S_H \setminus S_H\cap V)$.
\rm \item \it  The order of  
$\pi_{n-1}(B_H \setminus B_H \cap V) \otimes \bC$ 
divides $\prod_i 
\Delta_i$. 
\end{enumerate}
\end{theorem}
\begin{proof}
We may assume that $H$ is a member of a linear pencil $h : (\bC^{n+1}, 0) 
\to 
(\bC, 0)$ which scans 
the space $B$ and is in general position with respect to $V$. Denote $H_t = 
\{ h=t\}$.
(a). The space $S_{H_t} \setminus S_{H_t}\cap V$ is diffeomorphic to 
$S_{H_0} 
\setminus S_{H_0}\cap V$, due to the transversality of the sphere to $H_t 
\cap V$, for any $t$ within a small neighbourhood of $0$. This comes from 
the 
fact that the restriction $g_{|H_0}$ is an isolated singularity.

Now $S_{H_0} \setminus S_{H_0}\cap V$ is the total space of the Milnor 
fibration over the circle $S^1$, defined by the restriction of $g/\| g\|$.
 The claim follows from the homotopy exact sequence and from the fact that 
the Milnor fibre of the restriction $g_{|H_0}$ is homotopically a bouquet 
of 
spheres $\bigvee S^{n-1}$.
 
(c). We have the following commuting diagram:
%%%%%%%%%%%%%%%%
\begin{equation}\label{eq:dia}
 \begin{array}{ccc}
S_{H_0} \setminus S_{H_0}\cap V & \stackrel{\sim}{\hookrightarrow} & 
B_{H_0} 
\setminus 
B_{H_0}\cap V \\
\updownarrow & \ & \downarrow \\
S_{H_t} \setminus S_{H_t}\cap V & \hookrightarrow & B_{H_t} \setminus 
B_{H_t}\cap V 
\end{array}
\end{equation}
 The embedding on the first line is a homotopy equivalence, which comes 
from 
the local cone structure. The arrow at the left is a diffeomorphism, as 
shown 
above. The arrow to the right 
  is also an embedding. Moreover, we claim that the space $B_{H_t} 
\setminus 
B_{H_t}\cap V$ is homotopy equivalent to the space $B_{H_0} \setminus 
B_{H_0}\cap V$ to which one attaches cells of dimension $n$.
For example one can argue along the following lines using 
the map $\phi =(h,g) : (\bC^{n+1}, 0) \to \bC$ (cf. e.g. 
\cite{Le} and \cite{Ti-bri}).
 Let $\Gamma (h,g)$ be the polar curve of $g$ with respect to $h$ and let $\Delta$ denote the image of $\Gamma (h,g)$ by $\phi$,
see Figure \ref{f:1}.

\begin{figure}[hbtp]\label{f:polar}
\begin{center}
\epsfxsize=5cm
\leavevmode
\epsffile{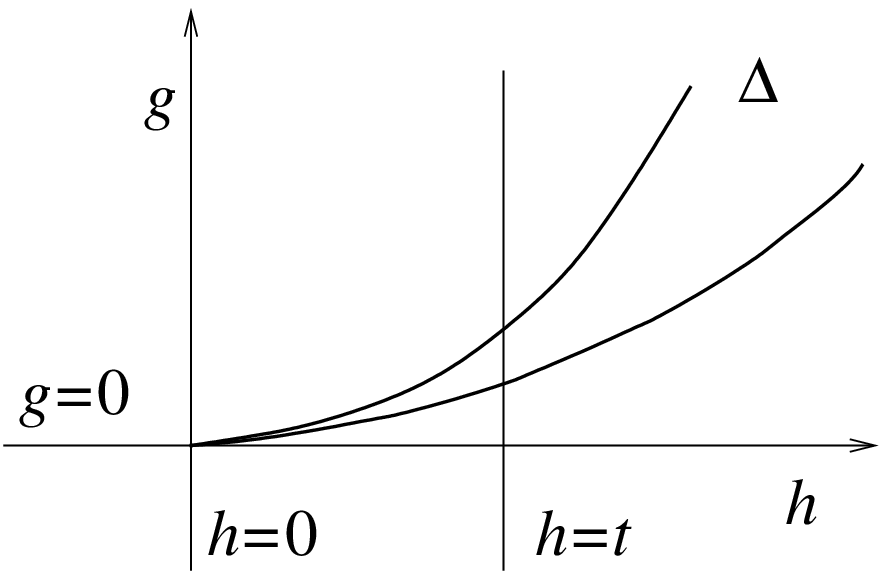}
\end{center}
\caption{The image by the map $(h,g) : (\bC^{n+1}, 0) \to (\bC^2, 0)$.} 
\label{f:1}   
\end{figure}

For suitable $\varepsilon$ and $t$, one has the following homotopy equivalences:
$B_{H_0} \setminus B_{H_0}\cap V \stackrel{\h}{\simeq} \{ |g| = \varepsilon\}\cap B_{H_0}$, by retraction, and  $\{ |g| = \varepsilon\}\cap B_{H_0}\stackrel{\h}{\simeq} \{ |g| = \varepsilon\}\cap B_{H_t}$, by isotopy.
The function $g$ has isolated singularities on $B_{H_t} 
\setminus B_{H_t}\cap V$, which are, by definition, exactly the points where the polar curve $\Gamma (h,g)$ cuts this space. It then follows that, up to homotopy equivalence, the space $B_{H_t} 
\setminus B_{H_t}\cap V$ can be constructed by attaching to  
$\{ |g| = \varepsilon\}\cap B_{H_t}$ a certain number of $n$-cells for each singular point of the function $g$ on $g^{-1}(D_\varepsilon) \cap B_{H_t} 
\setminus B_{H_t}\cap V$.

 Now, we apply $\pi_{n-1}$ to the diagram (\ref{eq:dia}). The arrow to the right  becomes a surjection in $\pi_{n-1}$ and this implies that the arrow on the 
bottom is also a surjection in $\pi_{n-1}$. By the multiplicativity of the 
order in exact sequences, it follows that $\Delta(B_H \setminus B_H \cap 
V)$ 
divides $\Delta(S_H \setminus S_H\cap V)$.
\end{proof} 

%%%%%%%%%%%%%%%%%%
%%%%%%%%%%%%%%%%
  The local analogue of Theorem \ref{t:alex} looks as follows:

\begin{theorem} 
\label{t:local}
Let $g$ be a germ of a holomorphic function with 
$\dim \Sing g \le 1$. Then the characteristic polynomial 
of the monodromy of $g$ acting on $H_{n-1}(M_g,\bC)$ 
coincides with the order of $\pi_{n-1}(B_H \setminus B_H \cap V)$.

\end{theorem}
\begin{proof} The proof follows from the Zariski-Lefschetz theorem 
of Hamm and L\^ e \cite{HL} and the equivariant identification of the universal cover of  
$B\setminus g^{-1}(0)$ with the Milnor fibre $M_g = g^{-1}(t)$, similar to the one used in the proof of \ref{t:alex}.
\end{proof}
 
\noindent 
{\it Acknowledgements.} The first named author thanks A. Dimca for pointing out an omission in a preliminary version 
of the statement of Corollary \ref{c:div} and for discussions related to this paper, in which he showed that he also has a proof of Proposition \ref{prop:gen}.

The second named author thanks Norbert A'Campo for suggesting 
the cup-product statement in Proposition \ref{prop:gen}.
%%%%%%%%%%%%%%%
%

\end{document}